\documentclass[a4paper,12pt]{article}
\usepackage[english]{babel}
\usepackage[T2A]{fontenc}
\usepackage[cp1251]{inputenc}
\usepackage{amsthm}
\usepackage[tbtags]{amsmath}
\usepackage{amsfonts,amssymb}
\begin{document}

\newtheorem{theorem}{Theorem}
\newtheorem{lemma}{Lemma}
\newtheorem{proposition}{Proposition}
\newtheorem{Cor}{Corollary}

\begin{center}
{\large\bf Centrally Essential Semirings}
\end{center}
\begin{center}
Lyubimtsev O.V., Tuganbaev A.A.
\end{center}


\textbf{Abstract.} A semiring is said to be centrally essential if for every non-zero element $x$, there exist two non-zero central elements $y, z$ with $xy = z$. We give some examples of non-commutative centrally essential semirings and describe some properties of additively cancellative centrally essential semirings.

The work of O.V. Lyubimtsev is done under the state assignment No~0729-2020-0055. A.A.Tuganbaev is supported by Russian Scientific Foundation, project 16-11-10013P.
																	
\section{Introduction}
By a \textsf{semiring}, we mean a structure that differs from an associative ring, possibly, by the irreversibility of the additive operation. In a semiring $S$, the zero is multiplicative by definition: we have $0s = s0 = 0$ for every $s\in S$. In our paper, we consider only semirings with $1$. For a semiring $S$, the \textsf{center} of $S$ is the set $C(S) = \{s\in S\colon ss' = s's$ for all $s'\in S\}$. This set is not empty, since it contains $0$ and $1$; we also have that $C(S)$ is a subsemiring in $S$. A semiring is said to be \textsf{centrally essential} if for every non-zero element $x$, there exist non-zero central elements $y, z$ with $xy = z$. 

Centrally essential rings with non-zero $1$ are studied, for example, in  \cite{MT18}, \cite{MT19}, \cite{MT19b}, \cite{MT20b}, and \cite{MT20c}. Every centrally essential semiprime ring with $1\ne 0$ is commutative; see \cite[Proposition 3.3]{MT18}. If $F$ is a field of order $2$ and $Q_8$ is the quaternion group of order 8, then the group algebra $FQ_8$ is a finite non-commutative centrally essential ring; see \cite{MT18}.
In \cite{MT20c}, a centrally essential ring $R$ is constructed such that the factor ring $R/J(R)$ with respect to the Jacobson radical is not a $PI$ ring (in particular, the ring $R/J(R)$ is not commutative). Matrix centrally essential algebras are studied in \cite{LT20la}. 
Abelian groups with centrally essential endomorphism rings are considered in \cite{LT20} and \cite{LT20lb}.

\textbf{1.1. Example.} We consider a semigroup $(M, \cdot)$ with multiplication table
\begin{center}
\begin{tabular}{ |c|c|c|c|c|c|}
 \hline
$\cdot$ & $0$ & $1$ & $a$ & $b$ & $c$  \\
\hline
$0$ & $0$ & $0$ & $0$ & $0$ & $0$ \\
\hline
$1$ & $0$ & $1$ & $a$ & $b$ & $c$ \\
\hline
$a$ & $0$ & $a$ & $a$ & $a$ & $c$ \\
\hline
$b$ & $0$ & $b$ & $b$ & $b$ & $c$ \\
\hline
$c$ & $0$ & $c$ & $c$ & $c$ & $c$ \\
 \hline
\end{tabular}
\end{center}

For a quick test of associativity, it is convenient to use the Light's associativity test; see \cite[p.7]{CP61}. Let $S = Sub\,(M)$ be the set of all subsets of the semigroup $M$.  For any $A, B\in S$, $S$, operations $A + B = A\cup B$ and $AB = \{ab\,|\, a\in A, b\in B\}$ are defined; then $S$ is a semiring with zero $\emptyset$ and the identity element $1 = 1_M$; see \cite[Example 1.10]{G99}. We have $|S| = 2^5 = 32$. We note that $S$ does not contain zero sums, i.e., the relation $A + B = \emptyset$ implies the relation $A = B = \emptyset$. In addition, $S$ is additively idempotent and multiplicatively idempotent. The center $C(S)$ os of the form
$$
C(S) = \{\emptyset, \{0\}, \{1\}, \{0, 1\}, \{c\}, \{0, c\}, \{0, 1, c\}, \{1, c\}\}.
$$
If $A\in S\backslash C(S)$, then $\emptyset\neq A\cdot \{c\}\in C(S)$. Consequently, $S$ is a non-commutative centrally essential semiring.

A semiring $S$ is said to be \textsf{reduced} if  $x = y$ for all $x, y\in S$ with $x^2 + y^2 = xy + yx$; see \cite{Ch93}. If $S$ is a ring, this is equivalent the property that $S$ has no non-zero nilpotent elements.
A semiring $S$ is said to be \textsf{additively cancellative} if the relation $x + z = y + z$ is equivalent to the relation $x = y$ for all $x, y, z\in S$. A ring $D(S)$ is called the \textsf{ring of differences} of the semiring $S$ if $S$ is a subsemiring in $D(S)$ and every element $a\in D(S)$ is the difference $x - y$ of some elements $x,y\in S$. The class of additively cancellative semirings contains all rings. The ring of differences is unique up to isomorphism over $S$; see \cite[Chapter II]{HW98} for details. An element $a$ of the semiring $R$ is called a \textsf{left zero-divisor} if $ab = 0$ for some $0\neq b\in S$. Similar to \cite[Lemma 2.2]{MT20b}, it can be proved that one-sided zero-divisors are two-sided zero-divisors in a centrally essential semiring. Other semiring-theoretical notions and designations can be found in \cite{G99, HW98}.

In the paper, we study properties of additively cancellative centrally essential semirings. The main result of the paper is the following theorem.

\textbf{1.2. Theorem.} There exists a non-commutative additively cancellative reduced centrally essential semiring without zero-divisors.
An additively cancellative reduced semiring $S$ is commutative if and only if the ring of differences of $S$ is a centrally essential ring.

\section{Additively Cancellative Centrally Essential Semirings}

A semiring $S$ is said to be \textsf{semiprime} if $S$ does not have nilpotent ideals. A semiring $S$ is said to be \textsf{semisubtractive} if for all $a, b\in S$ with $a\neq b$, there exists an element $x\in S$ such that $a + x = b$ or $b + x = a$.

\textbf{2.1. Proposition.} Let $S$ be an additively cancellative semisubtractive centrally essential semiring with center $C$. The following conditions are equivalent.
\begin{itemize}
\item
$S$ is a semiprime semiring.
\item
$C$ is a semiprime semiring.
\item
$S$ does not have non-zero nilpotent elements.
\item
$S$ is a commutative semiring without non-zero nilpotent elements.
\end{itemize}

\textbf{Proof.} It is well known that a semiring $S$ can be embedded in the ring of differences $D(S)$ if and only if $S$ is additively cancellative. In addition, the relation $D(S) = -S\cup S$ holds if and only if $S$ is a semisubtractive semiring; see \cite[Chapter II, Remark 5.12]{HW98}. Then the assertion follows from \cite[Proposition 2.8]{MT19b}.~\hfill$\square$
 
\textbf{2.2. Remark.} It follows from Example 1.1 that the assertion of Proposition 2.1 is not true without the assumptions of additive cancellativity and semisubtractivity.
In Example 3.2\footnote{See below.}, a non-commutative centrally essential semiring without zero-divisors is constructed; this semiring is additively cancellative but is not semisubtractive.

It is known that every idempotent of a centrally essential ring is central; see  \cite[Lemma 2.3]{MT18}. For semirings, a similar result is not true; see Example 1.1. For a semiring $S$, an idempotent $e$ of $S$ is said to be \textsf{complemented} if there exists an idempotent $f\in S$ with $e + f = 1$.

\textbf{2.3. Proposition.} In an additively cancellative centrally essential semiring $S$, any complemented idempotent is central. 

\textbf{Proof.} Let $e^2 = e$ and $e + f = 1$ for some $f\in S$. Since $S$ is an additively cancellative semiring, it follows from $e = e + fe$ that $fe = 0$. Similarly, we have $ef = 0$. Let $x\in S$ and $xe\neq 0$. Then $x = ex + fx$ and $xe = exe + fxe$. 

First, we assume that $fxe = 0$, i.e., $xe = exe$. Since $x = xe + xf$, we have $ex = exe + exf$. If $exf\neq 0$, then there exist $c, d\in C(S)$ 
with $(exf)c =d\neq 0$. Then 
$$
0\neq d = ed = de = (exfc)e = (exc)fe = 0;
$$
this is a contradiction. Consequently, $exf = 0$ and $ex = xe = exe$.

Now let $fxe\neq 0$. Then $0\neq (fxe)c = d$ for some non-zero elements $c, d\in C(S)$. In this case, 
$$
0\neq d = de = ed = ef(xec) = 0;
$$
this is a contradiction.~\hfill$\square$

\textbf{2.4. Corollary.} If $S$ is an additively cancellative semiring, then the semiring $M_n(S)$ of all matrices and the semiring $T_n(S)$ of all upper triangular matrices over $S$ are not centrally essential for $n\ge 2$.

\textbf{Proof.} For the identity matrices of the above semirings, we have $E = E_{11} + \ldots + E_{nn}$, where $E_{11}, \ldots, E_{nn}$ are matrix units. It follows from \cite[Example 4.19]{G99} that $M_n(S)$ is an additively cancellative semiring. The idempotents $E_{11}, \ldots, E_{nn}$ are non-central complemented idempotents. Consequently, the semirings $M_n(S)$ and $T_n(S)$ are not centrally essential.~\hfill$\square$

As it was mentioned above, additively cancellative semirings $S$ coincide with semirings $S$ which can be embedded in the rings of differences $D(S)$ whose elements are of the form $x - y$, where $x,y\in S$. 

\textbf{2.5. Example.} We consider the semiring $S$ generated by the matrices
$$
\left(\begin{matrix}
\alpha & a & b\\
0 & \alpha & c\\
0 & 0 & \alpha\\
\end{matrix}\right), \left(\begin{matrix}
0 & 0 & b\\
0 & 0 & 0\\
0 & 0 & 0\\
\end{matrix}\right), \left(\begin{matrix}
0 & 0 & 0\\
0 & 0 & 0\\
0 & 0 & 0\\
\end{matrix}\right), \left(\begin{matrix}
\alpha & 0 & 0\\
0 & \alpha & 0\\
0 & 0 & \alpha\\
\end{matrix}\right), 
$$
where $\alpha, a, b, c\in \mathbb{Z}^{+}$.
Let $A = (a_{ij})$ and $B = (b_{ij})$, where $a_{12} = b_{23} = a$, $b_{12} = a_{23} = c$, $a\neq c$, and the remaining components are equal to each other. 
Then $AB\neq BA$, i.e., $S$ is a non-commutative semiring. It is directly verified that the center $C(S)$ consists of matrices of the form 
$$
\left(\begin{matrix}
\alpha & 0 & b\\
0 & \alpha & 0\\
0 & 0 & \alpha\\
\end{matrix}\right),
$$
where $\alpha, b\in \mathbb{Z}^{+}\cup \{0\}$. Since $0\neq AD\in C(S)$, where $0\neq A\in S\backslash C(S)$, $0\neq D\in C(S)$ with $\alpha = 0$, we have that $S$ is a non-commutative centrally essential semiring. However, the ring of differences $D(S) = M_3(\mathbb{Z})$ is not a centrally essential ring, since the ring has non-central idempotents. In addition, any centrally essential subalgebra of a local triangular  $3\times 3$ matrix algebra is commutative; this is proved in \cite{LT20la}.

We give an example of a centrally essential ring $R$ which is the ring of differences for two proper semirings $S_1$ and $S_2$ of $R$ such that $S_1$ is not a centrally essential semiring and $S_2$ is a centrally essential semiring.

\textbf{2.6. Example.} Let $R$ be the ring consisting of matrices of the form 
$$
\left(\begin{matrix}
\alpha & a & b & c & d & e & f\\
0 & \alpha & 0 & b & 0 & 0 & d\\
0 & 0 & \alpha & 0 & 0 & 0 & e\\
0 & 0 & 0 & \alpha & 0 & 0 & 0\\
0 & 0 & 0 & 0 & \alpha & 0 & a\\
0 & 0 & 0 & 0 & 0 & \alpha & b\\
0 & 0 & 0 & 0 & 0 & 0 & \alpha\\
\end{matrix}\right)                          \eqno(1)
$$
over the ring $\mathbb{Z}$ of integers. In \cite{LT20la}, it is proved that $R$ is a non-commutative centrally essential ring.  
Let $S_1$ be the semiring generated by matrices of the form $(1)$ over $\mathbb{Z}^{+}$ and scalar matrices with 
$\alpha\in \mathbb{Z}^{+}\cup \{0\}$ and zeros на the remaining positions. 
Since $C(S_1)$ consists of scalar matrices, $S_1$ is not a centrally essential semiring. We note that $S_1$ is a semiring without zero-divisors. At the same time, the semiring $S_2$ of matrices of the form $(1)$ over the semiring $\mathbb{Z}^{+} \cup \{0\}$ is a centrally essential semiring.

\textbf{2.7. Lemma} \cite[Chapter II, Theorem 5.13]{HW98} In a semiring $S$, any central element is contaned in the center $C(D(S))$ of its ring of differences.

\textbf{2.8. Proposition.} Let $S$ be a centrally essential semiring without zero-divisors. If the ring $D(S)$ does not contain zero-divisors, the semiring $S$ is commutative.

\textbf{Proof.} Let $0\neq a = x - y\in D(S)$. By assumption, $0\neq xc = d$ and $0\neq yf = g$ for some $c, d, f, g\in C(S)$. Then
$$
a(cf) = (x - y)cf = (xc)f - (yf)c = df - gc.
$$
It follows from Lemma 2.7 that $c, d, f, g\in C(D(S))$ and $ac'\in C(D(S))$, where $c' = cf$. In addition, $ac'\neq 0$, since $D(S)$ does not contain zero-divisors. Then $D(S)$ is a commutative ring; see \cite[Proposition 3.3]{MT18}.~\hfill$\square$

\section{Proof of Theorem 1.2}

We recall that the \textsf{upper central series} of a group $G$ is the chain of subgroups 
$$
\{e\} = C_0(G)\subseteq C_1(G)\subseteq \ldots,
$$
where $C_i(G)/C_{i-1}(G)$ is the center of the group $G/C_{i-1}(G)$, $i\ge 1$. For a group $G$, the  \textsf{nilpotence class} of $G$ is the least positive integer $n$ with $C_n(G) = G$ provided such an integer $n$ exists.

\textbf{3.1. Proposition; cf. \cite[Proposition 2.6]{MT18}.} 
Let $G$ be a finite group of nilpotence class $n\le 2$ and let $S$ be a commutative semiring without zero-divisors or zero sums. Then $SG$ is a centrally essential group semiring.

\textbf{Proof.} If $n = 1$, then the group $G$ is Abelian and $SG$ is a centrally essential group semiring; see \cite[Lemma 2.2]{MT18}.

Let $n = 2$. Similar to the case od group rings (e.g., see \cite[Part 2]{P77}), the center $C(SG)$ is a free $S$-semimodule with basis
$$
\left\{\sum_{K} \, | \, K \mbox{are the conjugacy classes in the group $G$}\right\}.
$$
It is sufficient to verify that $SG\sum_{C(G)}\subseteq C(SG)$, where $C(G)$ is the center of the group $G$. Indeed, if $g, h\in G$, then 
$$
(gh)^{-1}hg\sum_{C(G)} = \sum_{C(G)},
$$
since $h^{-1}g^{-1}hg\in G'\subseteq C(G)$.~\hfill$\square$

We give an example of a noncommutative additively cancellative reduced centrally essential semiring without zero-divisors.

\textbf{3.2. Example.} Let $Q_8$ be the quaternion group, i.e., the group with two generators $a$, $b$ and defining relations 
$a^4 = 1$, $a^2 = b^2$ and $aba^{-1} = b^{-1}$; e.g., see \cite[Section 4.4]{H59}. We have 
$$
Q_8 = \{e, a, a^2, b, ab, a^3, a^2b, a^3b\},
$$
the conjugacy classes of $Q_8$ are
$$
K_e = \{e\}, K_{a^2} = \{a^2\}, K_a = \{a, a^3\}, K_b = \{b, a^2b\}, K_{ab} = \{ab, a^3b\},
$$
and the center $C(Q_8)$ is $\{e, a^2\}$. We consider the group semiring $SQ_8$, where $S = \mathbb{Q}^{+}\cup \{0\}$. 
Since $Q_8$ is a group of nilpotence class $2$, it follows from Proposition 3.1 that $SQ_8$ is a centrally essential group semiring.
To illustrate the above, we have 
$$
a\sum_{C(Q_8)} = \sum_{K_a},\quad b\sum_{C(Q_8)} = \sum_{K_b},
$$
$$
ab\sum_{C(Q_8)} = \sum_{K_{ab}},\quad a^3\sum_{C(Q_8)} = \sum_{K_a},
$$ 
$$
a^2b\sum_{C(Q_8)} = \sum_{K_b},\quad a^3b\sum_{C(Q_8)} = \sum_{K_{ab}}.
$$
The group ring of differences $\mathbb{Q}Q_8$ is a reduced ring; see \cite[Theorem 3.5]{S75}. Then $SQ_8$ is a reduced semiring. Indeed, if $x^2 + y^2 = xy + yx$ and $x\neq y$, then $x^2 + y^2 - xy - yx = (x - y)^2 = 0$ in the ring $\mathbb{Q}Q_8$; this is not true. Thus, $SQ_8$ is a non-commutative reduced centrally essential semiring without zero-divisors.  
We note that $\mathbb{Q}Q_8$ is not a centrally essential ring, since centrally essential reduced rings are commutative.~\hfill$\square$

\textbf{3.3. Completion of the Proof of Theorem 1.2.}

It follows from Example 3.2 that there exists a non-commutative additively cancellative reduced centrally essential semiring without zero-divisors.

If a semiring $S$ is commutative, then $D(S)$ is a commutative ring, i.e., $D(S)$ is centrally essential. Conversely, let $D(S)$ be a centrally essential ring. Since $S$ is a reduced semiring, $D(S)$ is a reduced ring. Indeed, let $0\neq a = x - y\in D(S)$.
If $a^2 = 0$, then $x^2 + y^2 = xy + yx$. Therefore, $x = y$, $a = 0$, and we have a contradiction. Then the ring $D(S)$ is commutative, since $D(S)$ is a reduced centrally essential ring. Consequently, $S$ is a commutative semiring.~\hfill$\square$

\section{Remarks and Open Questions}

An element $x$ of the semiring $S$ is said to be \textsf{left (resp., right) multiplicatively  cancellative} if $y = z$ for all $y, z\in S$ with $xy = xz$ (resp., $yx = zx$).
A semiring $S$ is said to be \textsf{left} (resp., \textsf{right})  \textsf{multiplicatively cancellative} if every $x\in  S\setminus \{0\}$ is left (resp., right) multiplicatively cancellative. A left and right  multiplicatively cancellative semiring is said to be \textsf{multiplicatively cancellative}; e.g., see \cite[Chapter I]{HW98}. 

\textbf{4.1. Remark.} A  left (resp., right) multiplicatively cancellative centrally essential semiring $S$ is commutative.

Indeed, let $a$ and $b$ be two non-zero elements of the semiring $S$. Since $S$ is a centrally essential semiring, there exists $c\in C(S)$ with $0\neq ac\in C(S)$. A left multiplicatively cancellative semiring does not contain left zero-divisors; see \cite[Chapter I, Theorem 4.4]{HW98}. Therefore, $acb\neq 0$. Then
$$
(ac)b = c(ab) = (ca)b = b(ca)= c(ba),
$$
whence we have $ab = ba$. A similar argument is true for right multiplicatively cancellative semirings.

A semiring with division, which is not a ring, is caled a  \textsf{division semiring}. A commutative division semiring is called a \textsf{semifield}.
It follows from Remark 4.1 that any centrally essential division semiring is a semifield. Indeed, it follows from \cite[Chapter I, Theorem 5.5]{HW98} that a division semiring with at lest two elements is multiplicatively cancellative.

\textbf{4.2. Open question.} Are there any non-commutative semisubtractive centrally essential semirings without non-zero nilpotent 
elements\footnote{See Proposition 2.1.}?

\textbf{4.3. Open question.} For groups of nilpotence class $n > 2$, are there any non-commutative centrally essential group semirings without zero-divisors?


\begin{thebibliography}{99}

\bibitem{Ch93} Chermnykh V.V. Sheaf representations of semirings// Russian Mathematical Surveys. -- 1993. -- Vol. 48, no. 5. -- P.~169-170.

\bibitem{CP61} Clifford A.H., Prieston G.B. The Algebraic Theory of Semigroups, vol. 1. AMS Survey No. 7, Providence, 1961.

\bibitem{G99} Golan J.S. Semirings and their applications. Kluwer Academic Publishers, Dordrecht; Boston; London, 1999.

\bibitem{H59} Hall M. The Theory of Groups, Macmillan, New York, 1959.

\bibitem{HW98} Hebisch U.,  Weinert H. J. Semirings. Algebraic theory and applications in computer science. World Scientific Publishing. Singapore, 1998.

\bibitem{LT20} Lyubimtsev O.V., Tuganbaev A.A. Centrally essential endomorphism rings of abelian groups// Comm. Algebra. --  2020. -- Vol. 48, no. 3. -- 
P.~1249-1256.

\bibitem{LT20la} Lyubimtsev O.V., Tuganbaev A.A. Local centrally essential subalgebras of triangular algebras// Linear and Multilinear Algebra. --
published on-line, https://doi.org/10.1080/03081087.2020.1802402.

\bibitem{LT20lb} Lyubimtsev O.V., Tuganbaev A.A. Centrally Essential Torsion-Free Rings of Finite Rank// Beitr\"age zur Algebra und Geometrie / Contributions to Algebra and Geometry. -- published on-line, https://doi.org/10.1007/s13366-020-00529-0.

\bibitem{MT18} Markov V.T., Tuganbaev A.A.  Centrally essential group algebras// J.~Algebra. -- 2018. -- Vol. 512, no. 15. -- P.~109-118.

\bibitem{MT19}  Markov V.T., Tuganbaev A.A. Centrally essential rings // Discrete Math. Appl. -- 2019. -- Vol. 29, no. 3. -- P.~189-194.

\bibitem{MT19b} Markov V.T., Tuganbaev A.A. Rings essential over their centers// Comm. Algebra. -- 2019. -- Vol. 47, no. 4. -- P.~1642-1649.

\bibitem{MT20b} Markov V.T., Tuganbaev A.A. Uniserial Noetherian Centrally Essential Rings// Comm. Algebra. -- 2020. -- Vol. 48, no. 1. -- P.~149-153.

\bibitem{MT20c} Markov V.T., Tuganbaev A.A. Constructions of Centrally Essential Rings// Comm. Algebra. -- 2020. -- Vol. 48, no. 1. -- P.~198-203.

\bibitem{P77} Passman D.S. The Algebraic Structure of Group Rings// John Wiley and Sons, New York, 1977.

\bibitem{S75} Sehgal K. S. Nilpotent elements in group rings// manuscripta mathematica. -- 1975. -- Vol. 15, no. 1 -- P.~65–80.

\end{thebibliography}
\end{document}